\documentclass[11pt]{article}
\usepackage{amsmath}
\usepackage{amsfonts}
\usepackage{amsthm}
\usepackage{verbatim}
\def\F{\mathbb{F}}

\begin{document}
\centerline{\bf{ A NOTE ON HYPERQUADRATIC CONTINUED FRACTIONS }}
\centerline{\bf{ IN CHARACTERISTIC 2 WITH PARTIAL QUOTIENTS OF DEGREE 1 }}
\vskip 0.5 cm
\centerline{\bf{by}}
\centerline{\bf{ A. Lasjaunias }}
\vskip 0.5 cm {\bf{Abstract.}} In this note, we describe a family of particular algebraic, and
nonquadratic, power series over an arbitrary finite field of characteristic 2, having a continued fraction expansion with all partial quotients of degree one. The main purpose being to point out a common origin with another analogue family in odd characteristic, which has been studied in [LY1].  
\vskip 0.5 cm
\noindent{\bf{1. Introduction }}

\par Let $p$ be a prime number, $q=p^s$ with $s\geq 1$, and let $\F_q$ be the finite field with $q$ elements. We let $\F_q[T]$, $\F_q(T)$
 and $\F(q)$  respectively denote,  the ring of polynomials, the field of rational functions and the field of power series in $1/T$ with coefficients in $\F_q$, where $T$ is a formal indeterminate. These fields are
valuated by the ultrametric absolute value (and its extension) introduced on $\F_q(T)$ by
$\vert P/Q\vert=\vert T\vert^{\deg(P)-\deg(Q)}$, where $\vert T\vert>1$ is a fixed real number. Hence a non-zero element of $\mathbb{F}(q)$ is written as $\alpha =\sum_{k\leq k_{0}}a_{k}T^{k}$ with $k_0\in \mathbb{Z}$, 
$a_{k}\in \mathbb{F}_{q}$, and $a_{k_{0}}\neq 0$ and we have $\vert \alpha\vert=\vert T\vert^{k_0}$. This field  $\F(q)$ is the completion of the field $\F_q(T)$ for this absolute value.

\par We recall that each
irrational (rational) element $\alpha$ of $\F(q)$ can be expanded as an infinite
(finite) continued fraction. This will be denoted
$\alpha=[a_1,a_2,\dots,a_n,\dots]$ where the $a_i\in \F_q[T]$, with
$\deg(a_i)>0$ for $i>1$, are the partial quotients. In this note we are concerned with infinite continued fractions in $\F(q)$ which are
 algebraic over $\F_q(T)$. Indeed we consider algebraic elements of a particular type. Let $r=p^t$ with 
$t\geq 0$, we say that $\alpha$ belonging to $\F(q)$ is hyperquadratic (of order $t$) if $\alpha$ is irrational and satisfies 
an algebraic equation of the particular form $A\alpha^{r+1}+B\alpha^r+C\alpha+D=0$, where $A,B,C$ and $D$ 
belong to $\F_q[T]$. For more details on these particular power series, the reader may see the introduction of [BL]. Note that quadratic elements are hyperquadratic. We recall that an element in $\F(q)$ is quadratic if and only if its sequence of partial quotients is ultimately periodic. Many explicit continued fractions are known for nonquadratic but hyperquadratic elements, see the references below. However, these particular algebraic elements play an important role in diophantine approximation and this is why they were first considered. For a general account on continued fractions in power series fields and diophantine approximation, also for more references on this matter, the
reader can consult [S] and [L2].

\par The origin of the study of continued fractions for hyperquadratic power series is due to Baum and Sweet ([BS1] and [BS2]). It has been developped in the 1980's by Mills and Robbins [MR]. Mills and Robbins pointed out the existence of such hyperquadratic continued fractions with all partial quotients of degree one, in odd characteristic with a prime base field. Different approaches were developped to explain this phenomenon (see [L1] and [LR1], [LR2]), first with the base field $\F_3$, then in all characteristic with an arbitrary base field. However these methods left many unanswered questions and the study could not be achieved. In [L3] a new approach was introduced and this led to the presentation of different families of explicit hyperquadratic continued fractions, also with unbounded partial quotients. In a recent work with Yao [LY1], with this new approach, in the case of an arbitrary base field of odd characteristic, we could give the description of a large family, including the historical examples due to Mills and Robbins, of hyperquadratic continued fractions with all partial quotients of degree one. 

\par As often, the case of even characteristic appears to be singular. The simplest case, where the base field is $\F_2$, has been fully treated by Baum and Sweet in [BS2]. They could describe the set of all the continued fractions with partial quotients of degree one (i.e. equal to $T$ or $T+1$). Among them, certain are algebraic and for every even integer $d\geq 2$ there is, in this set, an element algebraic of degree $d$. However, it was proved later that none of these algebraic elements (for $d>2$) were hyperquadratic (see [L2, p. 225]). Later, in a joint work with Ruch we could exhibit hyperquadratic continued fractions in $\F(2^s)$, with partial quotients all of the form $\lambda T$ where $\lambda \in \F_{2^s}$ (see [LR1, p. 280], with $r=2$ and also [LR2, p. 556]). Other examples were introduced by Thakur [T, p. 290]. Finally, we shall not omit the impressive note left by Robbins on arXiv [R], in which several conjectures on continued fractions for particular cubic power series in $\F(4)$ and $\F(8)$, with partial quotients of bounded degree, are presented.

\par Here we present a family of hyperquadratic (for each order $t\geq 1$) continued fractions in $\F(2^s)$ of the form $\alpha=[\lambda_1T,\lambda_2T,....,\lambda_nT,...]$. In the next section, we state a theorem showing how this sequence $(\lambda_n)_{n\geq 1}$ in $\F_{2^s}^*$ is defined recursively from an arbitrary number of initial values. Moreover, in a final remark, we make a comment on the poperties of this sequence. However, the main goal of this note is to show that this family has a common origin with another family of continued fractions in the case of odd characteristic, introduced earlier anf fully studied in [LY1]. Indeed, eventhough the case of characteristic $2$ is simpler, it appears that the existence of both families is directly connected to a particular universal quadratic power series $\omega=[T,T,...,T,....]$ belonging to $\F(p)$ for all prime numbers $p$. Note that this element $\omega$ is actually the analogue, in the formal case, of the celebrated quadratic real number $[1,1,\dots,1,\dots]=(1+\sqrt{5})/2$.

\vskip 0.5 cm
\noindent{\bf{2. Results }}

\par Our method, to obtain the explicit continued fraction for certain hyperquadratic power series, has been introduced in [L3]. We shall also use the notation, concerning continued fractions and continuants, presented in [LY1, Section 2].

\par Let $\alpha $ be an irrational element in $\mathbb{F}(q)$ with $\alpha
=[a_{1},\ldots ,a_{n},\ldots ]$ as its continued fraction expansion. We
denote by $\mathbb{F}(q)^{+}$ the subset of $\mathbb{F}(q)$ containing the
elements having an integral part of positive degree (i.e. with $\deg
(a_{1})>0$). For all integers $n\geq 1$, we put $\alpha
_{n}=[a_{n},a_{n+1},\ldots ]$ ($\alpha_1=\alpha$), called the $n$-th complete quotient. We introduce the usual 
continuants $x_{n},y_{n}\in \mathbb{F}_{q}[T]$ such that $x_{n}/y_{n}=[a_{1},a_{2},\ldots
,a_{n}] $ for $n\geq 1$. As usual we extend the latter notation to $n=0$ with $x_{0}=1$ and $y_{0}=0$. We recall that, with our notations for continuants, we will write $x_n=\langle a_1,a_2,\dots,a_n \rangle$ and $y_n=\langle a_2,\dots,a_n \rangle$. 

\par As above we set $r=p^{t}$, where $t\geq 0$ is an integer. Let $P,Q\in
\mathbb{F}_{q}[T]$ such that $\deg (Q)<\deg (P)<r$. Let $\ell \geq 1$ be an integer and 
$A_{\ell }=(a_{1},a_{2},\ldots ,a_{\ell })$ a vector in $(\mathbb{F}_{q}[T])^{\ell }$ such that $\deg (a_{i})>0$ for $1\leq
i\leq \ell $. Then by Theorem 1 in [L3, p. 333], there exists an infinite
continued fraction in $\mathbb{F}(q)$ defined by $\alpha
=[a_{1},a_{2},\ldots ,a_{\ell },\alpha _{\ell +1}]$ such that $\alpha
^{r}=P\alpha _{\ell +1}+Q$. This element $\alpha $ is hyperquadratic and it
is the unique root in $\mathbb{F}(q)^{+}$ of the following algebraic
equation:
$$y_{\ell }X^{r+1}-x_{\ell }X^{r}+(y_{\ell -1}P-y_{\ell }Q)X+x_{\ell
}Q-x_{\ell -1}P=0.  \eqno{(1)}$$
A continued fraction defined as above will be called of type $(r,\ell,P,Q)$. Note that such a continued fraction is depending on the choice of the vector $A_{\ell }=(a_{1},a_{2},\ldots ,a_{\ell })$ in $(\mathbb{F}_{q}[T])^{\ell }$.

\par Now, let us introduce an important sequence of polynomials in $\F_p[T]$, for all prime numbers $p$ (see [L3, p. 331]). This sequence $(F_n)_{n\geq 0}$ is defined by induction as follows:
$$F_0=1,\quad F_1=T \quad \text{ and}\quad F_{n+1}=TF_n+F_{n-1}\quad \text{ for} \quad n\geq 1.$$
Note the similarity with the celebrated sequence of Fibonacci integers. With our notation, for $n\geq 1$, we clearly have 
$$F_n/F_{n-1}=[T,T,\dots,T] \quad \text{ and}\quad F_n=\langle T,T,\dots,T \rangle ,$$ 
where the finite sequence of $T$ has length $n$. Then we define in $\F(p)$ the following infinite continued fraction 
$$\omega=[T,T,\dots,T,\dots]=\lim_{n\rightarrow \infty}F_n/F_{n-1}.$$
Note that $\omega$ is quadratic and satisfies $\omega=T+1/\omega$. We have the following lemma. 
\newline {\bf{Lemma 1. }}{\emph{ Let $p$ be a prime, the sequence $(F_n)_{n\geq 0}$ in $\F_p[T]$ and $\omega \in \F(p)$ be defined as above. As above, let $r=p^t$ where $t\geq 1$ is an integer. Then we have 
$$F_{n-1}=(\omega^n-(-\omega^{-1})^n)/(\omega +\omega^{-1})\quad \text{ for} \quad n\geq 1.$$  
$$\text{If }\quad p=2 \quad \text{then } \quad F_{r-1}=T^{r-1}\quad \text{and}\quad \text{ if }\quad p>2 \quad \text{then } \quad F_{r-1}=(T^2+4)^{(r-1)/2}.$$}}
\noindent Proof: For $n\geq 0$, we set $\Omega_n=\omega^n-(-\omega^{-1})^n$. Then we have
$$T\Omega_n+\Omega_{n-1}=\omega^{n-1}(T\omega+1)-(-\omega^{-1})^{n-1}(-T\omega^{-1}+1).\eqno{(2)}$$
Since $\omega=T+\omega^{-1}$, we get $\omega^2=T\omega+1$ and $\omega^{-2}=-T\omega^{-1}+1$. Hence $(2)$ becomes 
$$T\Omega_n+\Omega_{n-1}=\omega^{n+1}-(-\omega^{-1})^{n+1}=\Omega_{n+1}.\eqno{(3)}$$
Hence, by $(3)$, the sequence $(\Omega_{n})_{n\geq 0}$ satisfies the same recurrence as the sequence $(F_n)_{n\geq 0}$. We observe that $\Omega_1/(\omega+\omega^{-1})=1=F_0$ and $\Omega_2/(\omega+\omega^{-1})=\omega-\omega^{-1}=T=F_1$. Since the sequences $(F_n)_{n\geq 0}$ and $(\Omega_{n+1}/(\omega+\omega^{-1}))_{n\geq 0}$ satisfy the same linear recurrence relation and are equal on the first two values, they coincide and we have $F_{n-1}=\Omega_n/(\omega+\omega^{-1})$, for $n\geq 1$, as stated in the lemma. If $p=2$, using the Frobenius isomorphism, we have $\Omega_r=(\omega+\omega^{-1})^r$ and also $\omega+\omega^{-1}=T$. Therefore, we can write $F_{r-1}= (\omega+\omega^{-1})^r/(\omega+\omega^{-1})=(\omega+\omega^{-1})^{r-1}=T^{r-1}$. If $p>2$, again using the Frobenius isomorphism, we also have $\Omega_r=(\omega+\omega^{-1})^r$. But here $\omega-\omega^{-1}=T$ and $(\omega-\omega^{-1})^2+4=(\omega+\omega^{-1})^2$. Consequently, we get $F_{r-1}=(\omega+\omega^{-1})^{r-1}=((\omega-\omega^{-1})^2+4)^{(r-1)/2}=(T^2+4)^{(r-1)/2}$. So the proof is complete.

\par We are now interested in the family of continued fractions of type $(r,\ell,\epsilon_1F_{r-1},\epsilon_2F_{r-2})$, for a pair $(\epsilon_1,\epsilon_2)\in  (\F_q^*)^2$. Actually this family contains continued fractions with all partial quotients of degree one, hence the $\ell$ first partial quotients must be chosen as polynomial of degree 1 in $\F_q[T]$. Note that, in the first works on hyperquadratic continued fractions with bounded partial quotients, the role played by the pair $(F_{r-1},F_{r-2})$ was not put in evidence. There are two distinct cases according to the characteristic $p=2$ or $p>2$. The case of characteristic $p>2$ was the first to be considered. In [MR], Mills and Robbins described, in the case $q=r=p>3$ , with $\ell=2$ and $(a_1,a_2)=(aT,bT)$, a family of such examples (see also the introduction of [L3, p. 332]). However in that case, in general, the solution of $(1)$ will not have all partial quotients of degree 1. This will happen if a mysterious condition, between the coefficients of the polynomials in $A_{\ell}$ and the pair $(\epsilon_1,\epsilon_2)$, is satisfied. It has been long and complicated to obtain this condition (see [LY1] and particularly the first section for the background information on the matter). Here, we shall now only consider the case of characteristic 2. Due to the simpler form of $F_{r-1}$, given in the previous lemma, we will avoid the sophistication appearing in odd characteristic.
 \par The following general lemma, on the continued fraction algorithm, is the basic tool of our method. It can be found in several articles (see [L3, Lemma 3.1, p.336)] or [LY1, p.267]).
\newline{\bf{Lemma 2. }}{\emph{For $n\geq 2$, given $n+1$ variables $x_1,x_2,\dots,x_n$ and $x$, we have formally 
$$[[x_1,x_2,\dots,x_n],x]=[x_1,x_2,\dots,x_n,y],$$ where 
$$y=(-1)^{n-1}\langle x_2,\dots,x_n\rangle^{-2}x-\langle x_2,\dots,x_{n-1}\rangle\langle x_2,\dots,x_n\rangle^{-1}.$$}}
\par As a consequence of Lemma 1 and Lemma 2, we have now the following lemma.
\newline{\bf{Lemma 3. }}{\emph{Let $q=2^s$ and $r=2^t$ where $s,t \geq 1$ are integers. Let $\alpha =[a_1,a_2,\dots,a_n,\dots]$ be an infinite continued fraction in $\F(q)$. We assume that there is a pair $(\epsilon_1,\epsilon_2)\in  (\F_q^*)^2$ and a pair of integers $(i,j)$, where $1\leq i<j$, such that we have $\alpha_i^r=\epsilon_1F_{r-1}\alpha_j+\epsilon_2F_{r-2}$. If $a_i=\lambda_iT$ where $\lambda_i\in \F_q^*$, then we have 
$$a_j=\epsilon_1^{-1}\lambda_i^rT \quad \text{and }\quad a_{j+k}=(\epsilon_2/\epsilon_1)^{(-1)^k}T\quad \text{for }\quad 1\leq k\leq r-1$$
and also
$$\alpha_{i+1}^r=\epsilon_1\epsilon_2^{-2}F_{r-1}\alpha_{j+r}+\epsilon_2^{-1}F_{r-2}.$$}}
\noindent Proof: From the equality $\alpha_i^r=\epsilon_1F_{r-1}\alpha_j+\epsilon_2F_{r-2}$, and since $\alpha_i^r=[a_i,\alpha_{i+1}]^r=[a_i^r,\alpha_{i+1}^r]$, we obtain
$$[\frac{a_i^r+\epsilon_2F_{r-2}}{\epsilon_1F_{r-1}}, \epsilon_1F_{r-1}\alpha_{i+1}^r]=\alpha_j.\eqno{(4)}$$
We denote by $W_i$ the word $T,T,...T$, where $T$ is repeated $i$ times. We have $(\epsilon_1/\epsilon_2)F_{r-1}/F_{r-2}=(\epsilon_1/\epsilon_2)[W_{r-1}]$. By Lemma 1, we have $F_{r-1}=T^{r-1}$, therefore $(4)$ becomes 
$$[\epsilon_1^{-1}\lambda_i^rT,(\epsilon_1/\epsilon_2)[W_{r-1}], \epsilon_1F_{r-1}\alpha_{i+1}^r]=[[\epsilon_1^{-1}\lambda_i^rT,x_2,x_3,\dots,x_r],\epsilon_1F_{r-1}\alpha_{i+1}^r]=\alpha_j,\eqno{(5)}$$
where $x_i=(\epsilon_1/\epsilon_2)^{(-1)^i}T$ for $2\leq i\leq r$. We shall now apply Lemma 2. Using classical properties of continuants (see [LY1, p. 266]), we can write
$$\langle x_2,\dots,x_r\rangle=\langle (\epsilon_1/\epsilon_2)T,(\epsilon_2/\epsilon_1)T,\dots,(\epsilon_1/\epsilon_2)T\rangle =(\epsilon_1/\epsilon_2)\langle T,T,\dots,T\rangle=(\epsilon_1/\epsilon_2)F_{r-1}$$
and also 
$$\langle x_2,\dots,x_{r-1}\rangle=\langle (\epsilon_1/\epsilon_2)T,(\epsilon_2/\epsilon_1)T,\dots,(\epsilon_2/\epsilon_1)T\rangle =\langle T,T,\dots,T\rangle=F_{r-2}.$$
We set $x_1=\epsilon_1^{-1}\lambda_i^rT$ and $x=\epsilon_1F_{r-1}\alpha_{i+1}^r$. Thus, according to Lemma 2, we get
$$[[x_1,x_2,x_3,\dots,x_r],x]=[x_1,x_2,\dots,x_r,((\epsilon_1/\epsilon_2)F_{r-1})^{-2}x+((\epsilon_1/\epsilon_2)F_{r-1})^{-1}F_{r-2}].\eqno{(6)}$$
Combining $(5)$ and $(6)$, since $x=\epsilon_1F_{r-1}\alpha_{i+1}^r$, we obtain
$$\alpha_j=[x_1,x_2,\dots,x_r,(\epsilon_2^2/\epsilon_1)F_{r-1}^{-1}\alpha_{i+1}^r+(\epsilon_2/\epsilon_1)F_{r-1}^{-1}F_{r-2}]=[x_1,x_2,\dots,x_r,y].\eqno{(7)}$$
Since $\vert \alpha_{i+1}\vert >1$, we observe that $\vert y\vert >1$. Comparing $(7)$ to $\alpha_j=[a_j,a_{j+1},\dots,a_{j+r-1},\alpha_{j+r}]$, we get $a_j=x_1$, $a_{j+1}=x_2,\dots a_{j+r-1}=x_r$ and $\alpha_{j+r}=y$. With the values given to $x_i$, for $1\leq i\leq r-1$, we obtain the values for the partial quotients from $a_j$ to $a_{j+r-1}$. Finally, with the value of $y$ in $(7)$, we obtain the last formula of the lemma. So the proof is complete.
\par We can now state and prove the following theorem.
\newline{\bf{Theorem. }}{\emph{Let $q=2^s$ and $r=2^t$ where $s,t \geq 1$ are integers. Let $\ell \geq 1$ be an integer and the pair $(F_{r-1},F_{r-2})$ in $(\F_2[T])^2$ be defined as above. Let $\Lambda_{\ell+2}=(\lambda_1,\lambda_2,\dots,\lambda_{\ell},\epsilon_1,\epsilon_2)$ be given in $(\F_q^*)^{\ell+2}$. Let $x_\ell,y_{\ell},x_{\ell-1}$ and $y_{\ell-1}$ be the four continuants in $\F_q[T]$, built from the vector $A_{\ell}=(\lambda_1T,\lambda_2T,\dots,\lambda_{\ell}T)$. We consider the following algebraic equation
$$y_{\ell }X^{r+1}+x_{\ell }X^{r}+(\epsilon_1y_{\ell -1}F_{r-1}+\epsilon_2y_{\ell }F_{r-2})X+\epsilon_2x_{\ell }F_{r-2}
+\epsilon_1x_{\ell -1}F_{r-1}=0. \eqno{(E)}$$
Then $(E)$ has a unique root $\alpha$ in $\F(q)^{+}$ and we have $\alpha=[\lambda_1T,\lambda_2T,\dots,\lambda_{\ell}T,\dots,\lambda_nT,\dots]$. The sequence $(\lambda_n)_{n\geq 1}$ in $\F_q^*$ is defined recursively, from the vector $\Lambda_{\ell+2}$, as follows :
$$\lambda_{\ell+rm+1}=(\epsilon_2/\epsilon_1)\epsilon_2^{(-1)^{m+1}}\lambda_{m+1}^r\quad \text{and }\quad \lambda_{\ell+rm+i}=(\epsilon_1/\epsilon_2)^{(-1)^i}\quad \text{for }\quad m\geq 0 \quad \text{and }\quad 2\leq i\leq r.$$
}} 
\noindent Proof: We observe that $(E)$ is simply equation $(1)$, stated above, built from the vector $A_{\ell}=(\lambda_1T,\lambda_2T,\dots,\lambda_{\ell}T)$ and the pair $(P,Q)=(\epsilon_1F_{r-1},\epsilon_2F_{r-2})$. According to the theorem from [L3], mentioned before, we know that $(E)$ has a unique root $\alpha$ in  $\F(q)^{+}$. This element is the infinite continued fraction defined by $\alpha=[\lambda_1T,\lambda_2T,\dots,\lambda_{\ell}T,\alpha_{\ell+1}]$ and $\alpha^r=\epsilon_1F_{r-1}\alpha_{\ell+1}+\epsilon_2F_{r-2}$. Since $a_1=\lambda_1T$, we can apply Lemma 3, with $i=1$ and $j=\ell+1$. Hence, we get
$$a_{\ell+1}=\epsilon_1^{-1}\lambda_1^rT \quad \text{and }\quad a_{\ell+i}=(\epsilon_1/\epsilon_2)^{(-1)^i}T\quad \text{for }\quad 2\leq i\leq r.\eqno{(8)}$$
We also have
$$\alpha_{2}^r=\epsilon_1\epsilon_2^{-2}F_{r-1}\alpha_{\ell+r+1}+\epsilon_2^{-1}F_{r-2}.$$
Let us consider the application $f$ from $(\F_q^*)^2$ into $(\F_q^*)^2$, defined by  $f(x,y)=(xy^{-2},y^{-1})$. We observe that $f$ is an involution. Hence, we can apply again Lemma 3 (note that $a_2=\lambda_2T$ even if $\ell=1$), replacing $(i,j)$ by $(2,l+r+1)$ and the pair $(\epsilon_1,\epsilon_2)$ by the pair $f(\epsilon_1,\epsilon_2)$. We obtain
$$a_{\ell+r+1}=\epsilon_2^2\epsilon_1^{-1}\lambda_2^rT \quad \text{and }\quad a_{\ell+r+i}=(\epsilon_1/\epsilon_2)^{(-1)^i}T\quad \text{for }\quad 2\leq i\leq r\eqno{(9)}$$
and 
$$\alpha_{3}^r=\epsilon_1F_{r-1}\alpha_{\ell+r+1}+\epsilon_2F_{r-2}.$$
 Hence, by repeated application of this lemma, we see that $a_i=\lambda_iT$ for $i\geq 1$. Note that $(\epsilon_2/\epsilon_1)\epsilon_2^{(-1)^{m+1}}$ is either $\epsilon_1^{-1}$ or $\epsilon_2^2\epsilon_1^{-1}$, according to the parity of the integer $m\geq 0$. From $(8)$ and $(9)$, we get the formulas defining the sequence $(\lambda_n)_{n\geq 1}$ in $\F_q^*$, stated in the theorem. So the proof is complete.
\vskip 0.5 cm
\noindent{\bf{Remark. }} In a recent joined work with J-Y. Yao [LY2], we have observed that the sequence of the leading coefficients of the partial quotients for an hyperquadratic continued fraction could be an automatic sequence. Naturally, the question arises for the sequence described in this theorem. In the particular case $r=2$, with the terminology intoduced above, we note that the element $\alpha$ is of type $(2,\ell,\epsilon_1T,\epsilon_2)$. Indeed, we have proved in [LY2, Theorem 3] that the corresponding sequence is 2-automatic. Let us indicate another criterion for a sequence in $\F_q$ to be automatic (see [AS, p. 356]) : $(\lambda_n)_{n\geq 1}$ in $\F_q$ is p-automatic if the power series $\sum_{i\geq 1}\lambda_iT^{-i}\in \F(q)$ is algebraic over $\F_q(T)$. Concerning the sequence $(\lambda_n)_{n\geq 1}$ described in our theorem, we conjecture the following : Set $\theta=\sum_{i\geq 1}\lambda_iT^{-i}\in \F(q)$, then there exist two elements $A$ and $B$ in $\F_q(T)$, depending on $r$ and on the vector $\Lambda_{\ell+2}$, such that we have $\theta^r+A\theta+B=0$.
Accordingly, by the criterion cited, the conjecture would imply the 2-automaticity of the sequence $(\lambda_n)_{n\geq 1}$ for all $r=2^t$ and $t\geq 1$. This conjecture is based on a partial study of the sequence $(\lambda_n)_{n\geq 1}$. Indeed, in the simplest case $r=2$, the conjecture is true and this study was complete enough to allow us to obtain the following :
$$\theta=\sum_{i=1}^{\ell}\lambda_iT^{-i}+(\epsilon_1/\epsilon_2)T^{-\ell}(T+1)^{-2}(\epsilon_2^{(-1)^\ell}T+1)+(\epsilon_1/\epsilon_2)\epsilon_2^{(-1)^\ell}T^{\ell-1}\rho$$
and
$$\rho^2+\rho+T^{-2\ell}\sum_{i=1}^{\ell}(1+\lambda_i^2(\epsilon_2/\epsilon_1)^2\epsilon_2^{(-1)^i-(-1)^\ell})T^{-2i+2}=0.$$

\begin{tabular}{ll}
 Alain LASJAUNIAS\\
Institut de Math\'ematiques de Bordeaux  \\
CNRS-UMR 5251\\
Talence 33405 \\
France \\
E-mail: Alain.Lasjaunias@math.u-bordeaux1.fr\\

\end{tabular}%

\end{document}